\documentclass[a4paper,11pt,reqno]{amsart}
\usepackage{xcolor}
\usepackage{amsmath, amsfonts, amsthm, graphicx, hyperref}
\usepackage[T1]{fontenc}
\usepackage{lmodern}
\DeclareMathOperator{\Var}{Var}
\DeclareMathOperator{\Out}{Out}
\DeclareMathOperator{\Inn}{Inn}
\DeclareMathOperator{\Aut}{Aut}
\newcommand{\Z}{\mathbb{Z}}
\newcommand{\N}{\mathbb{N}}
\newcommand{\cA}{\mathcal{A}}
\newcommand{\cB}{\mathcal{B}}
\newcommand{\cC}{\mathcal{C}}
\newcommand{\cD}{\mathcal{D}}
\newcommand{\Luczak}{\L uczak}
\newcommand{\Swiatkowski}{\'Swi\k{a}tkowski}
\newcommand{\Odrzygozdz}{Odrzyg\'o\'zd\'z}

\numberwithin{equation}{section}
\newtheorem{theorem}[equation]{Theorem}
\newtheorem{proposition}[equation]{Proposition}
\newtheorem{corollary}[equation]{Corollary}
\newtheorem{lemma}[equation]{Lemma}

\theoremstyle{definition}
\newtheorem{definition}[equation]{Definition}

\title{Property FA for random $\ell$-gonal groups}
\author{Emily Clement}
\address{CNRS, LIPN UMR 7030, Université Sorbonne Paris Nord, F-93430 Villetaneuse, France, 
 \url{https://perso.eleves.ens-rennes.fr/people/emily.clement/} }
\email{emily.clement@lipn.univ-paris13.fr}
\author{John M. Mackay}
\address{School of Mathematics, University of Bristol, Bristol, BS8 1TX.}
\email{john.mackay@bristol.ac.uk}
\date{\today}
\begin{document}

\begin{abstract}
	In the binomial $\ell$-gonal model for random groups, where the random relations all have fixed length $\ell\geq 3$ and the number of generators goes to infinity, we establish a double threshold near density $d=\frac{1}{\ell}$, where the group goes from being free to having Serre's property FA.
	As a consequence, random $\ell$-gonal groups at densities $\frac{1}{\ell} < d< \frac{1}{2}$ have boundaries homeomorphic to the Menger sponge, and $\frac{1}{\ell}$ is also the threshold for finiteness of $\Out(G)$.  We also see that the thresholds for property FA and Kazhdan's property (T) differ when $\ell \geq 4$. 

	Our methods are inspired by work of Antoniuk--\Luczak--\Swiatkowski\ and Dahmani--Guirardel--Przytycki.
\end{abstract}

\maketitle

\section{Introduction}
In the study of random groups, as pioneered by Gromov~\cite{Gro93}, one considers the properties of randomly chosen finitely presented groups $\langle S | R \rangle$.  As with other random structures, finding `thresholds' for properties to hold is of particular interest.  In this paper, we study random $l$-gonal groups, where all relators have length $\ell$, and find a new double threshold for Serre's property FA.  These groups are typically Gromov hyperbolic and so have a boundary at infinity, and our results show that at this same threshold the boundary goes from being a Cantor set to being a Menger sponge, with an intermediate range where it is neither.  We also find the threshold for the finiteness of the random group's outer automorphism group.

\subsection{Main results}
The model of random group we consider is the following.

\begin{definition}[{cf.\ \cite{AntLucSwi-13threshold}, \cite{AshcroftRonDou20}}]
	\label{def-l-gonal-density}
	For each length $\ell \in \{3,4,\ldots\}$, probability $p \in [0,1]$ and number of generators $m \geq 2$, 
	let $\mathcal{M}(m,\ell,p)$ be the probability space consisting of group presentations $G = \langle S | R \rangle$ where $S = \{s_1,\ldots, s_m\}$, and $R$ consists of relators where each of the cyclically reduced words of length $\ell$ in $S^\pm$ are chosen independently with probability $p$.  
	We call a group (presentation) sampled from $\mathcal{M}(m,\ell,p)$ a \emph{(binomial) random $\ell$-gonal group with probability $p$ (and $m \geq 2$ generators)}.

	For a function $p=p(m)$, a random $\ell$-gonal group satisfies a given property \emph{asymptotically almost surely} (a.a.s.) if the probability $G = \langle S|R \rangle \in \mathcal{M}(m,l,p)$ satisfies the property goes to $1$ as $m \to \infty$.
\end{definition}
% I think maybe skip this sentence here?  Note that the number of cyclically reduced words of length $\ell$ in $m$ generators is $(1+o(1))(2m-1)^\ell$. 
If $p(m)= m^{\ell (d-1)}$ for some fixed \emph{density} $d \in (0,1)$, then this is a binomial variation on the $\ell$-gonal (or ``$\ell$-angular'') model of Ashcroft--Roney-Dougal~\cite{AshcroftRonDou20}; the triangular binomial model was studied by Antoniuk--\Luczak--\Swiatkowski~\cite{AntLucSwi-13threshold}.  By standard monotonicity results for random structures, our results will hold with a threshold for this choice of $p$ if and only if there is a threshold in the $\ell$-gonal models where the number of relations is fixed at $\lfloor m^{\ell d}\rfloor$ (see Section~\ref{sec-not-free-not-fa}).
Earlier studies of these related models were made by \.Zuk~\cite{Zuk03} ($\ell=3$) and \Odrzygozdz~\cite{Odr16, Odr19} ($\ell=4,6$), motivated in part by studying Kazhdan's Property (T).  

Random groups at low densities are more flexible than those at high densities.  For example, if $d<\frac{1}{\ell}$ a random $\ell$-gonal group is a.a.s.\ free, if $\frac{1}{\ell}<d<\frac{1}{2}$ it is a.a.s.\ infinite non-free torsion-free hyperbolic, and if $\frac{1}{2}<d$ it is trivial or $\Z/2\Z$~\cite{AntLucSwi-13threshold,AshcroftRonDou20}.
In the Gromov density model (where the length of relators goes to infinity), Dahmani--Guirardel--Przytycki show that at all densities $d>0$, random groups have a.a.s.\ Serre's property FA~\cite{DahGuiPrz11-nosplit}.
\begin{definition}
	\label{def-property-fa}
	A group $G$ has \emph{property FA} if every action of $G$ on a simplicial tree has a fixed point.
\end{definition}
A group has property FA if and only if it does not split non-trivially as an amalgamated product or HNN extension \cite{SerreTrees03}.
Groups in the $\ell$-gonal models are `less quotiented' than in the Gromov density model, so it is natural to wonder whether groups in the $\ell$-gonal model also have property FA.
Free groups do not have property FA, so neither do random $\ell$-gonal groups at density $d<\frac{1}{\ell}$.  Our main result shows that this is sharp, in fact we establish a double threshold for property FA near density $d=\frac{1}{\ell}$.
\begin{theorem}
	\label{thm-propFA-sharp}
	Fix $\ell \geq 3$.  
	There exists $c,C,c',C'>0$ so that for a random $\ell$-gonal group $G \in \mathcal{M}(m,l,p)$,
	\begin{enumerate}
		\item if $p \leq cm^{1-\ell}$ then a.a.s.\ $G$ is a non-abelian free group;
		\item if $Cm^{1-\ell} \leq p \leq c' \log(m) m^{1-\ell}$ then a.a.s.\ $G$ splits as $G = H \ast F_q$ for some $q \geq \frac{\sqrt{m}}{2}$ with $H$ not free, so in particular $G$ is not free and does not have property FA; and
		\item if $p \geq C' \log(m)m^{1-\ell}$ then a.a.s.\ $G$ satisfies Property FA.
	\end{enumerate}	
\end{theorem}

Any infinite Gromov hyperbolic group $G$ has a boundary at infinity $\partial_\infty G$, a canonically defined compact topological space.  Dahmani--Guirardel--Prz\-yty\-cki use property FA to show that the boundary of a random group in Gromov's density model is a.a.s.\ homeomorphic to a Menger sponge at densities $0<d<\frac{1}{2}$.  Likewise, we are able to deduce the following corollary.
\begin{corollary}
	\label{cor-boundary-sharp}
	Fix $\ell \geq 3$.  
	There exists $c,C,c',C'>0$ so that for a random $\ell$-gonal group $G \in \mathcal{M}(m,l,p)$,
	\begin{enumerate}
		\item if $p \leq cm^{1-\ell}$ then a.a.s.\ $\partial_\infty G$ is homeomorphic to a Cantor set;
		\item if $Cm^{1-\ell} \leq p \leq c' \log(m) m^{1-\ell}$ then a.a.s.\ $\partial_\infty G$ has infinitely many connected components but is not homeomorphic to a Cantor set;
		\item if for some $\epsilon>0$ we have $C' \log(m)m^{1-\ell}< p < m^{-\ell/2 - \epsilon}$ then a.a.s.\ $\partial_\infty G$ is homeomorphic to a Menger sponge.
	\end{enumerate}	
\end{corollary}
Of course, if $p > m^{-\ell/2+\epsilon}$ then a.a.s.\ $G$ is finite, so $\partial_\infty G = \emptyset$.
It is natural to wonder whether in case (2) the topological type is determined: one possibility is that $\partial_\infty G$ is a.a.s.\ homeomorphic to a ``Cantor set of points and Menger sponges'', i.e.\ to $\partial_\infty (H \ast \Z)$ where $H$ is a hyperbolic group with Menger sponge boundary.

Theorem~\ref{thm-propFA-sharp} has further consequences for the algebraic structure of random $\ell$-gonal groups, as in \cite[Corollary 1.6]{DahGuiPrz11-nosplit}.  Recall that for a group $G$, its \emph{outer automorphism group} $\Out(G)$ is the group of automorphisms of $G$, quotiented by the normal subgroup of inner automorphisms.  For any groups $G$ and $\Gamma$, $\mathrm{Hom}(G; \Gamma)$ is the set of homomorphisms $G \to \Gamma$; there is a natural equivalence relation on this set given by post-composition by conjugation in $\Gamma$.
\begin{corollary}
	\label{cor-finite-out-sharp}
	Fix $\ell \geq 3$.  
	There exists $c',C'>0$ so that for a random $\ell$-gonal group $G \in \mathcal{M}(m,l,p)$,
	\begin{enumerate}
		\item if $p \leq c' \log(m) m^{1-\ell}$ then a.a.s.\ $\Out(G)$ is infinite and for any infinite torsion-free hyperbolic group $\Gamma$, $\mathrm{Hom}(G; \Gamma)$ is not finite even up to conjugacy.
		\item if for some $\epsilon>0$ we have $C' \log(m)m^{1-\ell}< p < m^{-\ell/2 - \epsilon}$ then a.a.s.\ $\Out(G)$ is finite and for any torsion-free hyperbolic group $\Gamma$, $\mathrm{Hom}(G;\Gamma)$ is finite up to conjugacy.
	\end{enumerate}	
\end{corollary}

\subsection{Further context and remarks}
Theorem~\ref{thm-propFA-sharp} is essentially already known for the case $\ell=3$ (``random triangular groups'') by work of Antoniuk, \Luczak\ and \Swiatkowski~\cite[Theorems 1--3]{AntLucSwi-13threshold}.
Indeed our proof of Theorem~\ref{thm-propFA-sharp}(1,2) largely follows their argument.
However, for Theorem~\ref{thm-propFA-sharp}(3), we need a new approach.
Antoniuk--\Luczak--\Swiatkowski\ show that random triangular groups have Kazhdan's Property (T) for $p \geq C'\log(m)m^{-2}$.
Kazhdan's Property (T) is a strong rigidity property of a group, restricting its unitary representations.  (We do not require the definition of Property (T) in this paper; \cite{BekkaHV09} is the standard reference.)  Any group with Property (T) also has property FA, so \cite[Theorem 3]{AntLucSwi-13threshold} is stronger than Theorem~\ref{thm-propFA-sharp}(3) for $\ell=3$.
However, not every group with Property FA has Property (T): for example, hyperbolic triangle groups have FA but do not have Property (T) (since they have FA but not all finite index subgroups have FA, see \cite[6.3.5]{SerreTrees03}).  The threshold for which $\ell$-gonal random groups have Property (T) is not known in most cases, but bounds are known.
For the triangular ($\ell=3$) model, prior to \cite{AntLucSwi-13threshold}, \.Zuk~\cite{Zuk03} and Kotowski--Kotowski~\cite{Kotowski13} had shown that Property (T) held a.a.s.\ at densities strictly greater than $ \frac{1}{3}$ (and failed at densities strictly smaller than $\frac{1}{3}$).
\Odrzygozdz\ showed the same holds in the hexagonal ($\ell=6$) model~\cite{Odr19}.
Ashcroft showed that for $\ell = 3q-r \geq 3$, $r=0,1$ or $2$, a random $\ell$-gonal group has Property (T) for densities $d> \frac{1}{3}+\frac{r}{3\ell}$~\cite{Ashcroft23-prop-T-lgonal}.
So for these ranges of densities, Theorem~\ref{thm-propFA-sharp}(3) is known, but in all other cases it is new as Property (T) is not known to hold.

In fact, for every $\ell \geq 4$, $\ell \neq 5$, Theorem~\ref{thm-propFA-sharp}(3) applies for a range of densities where Property (T) is known \emph{not} to hold.  
For example, in the square ($\ell=4$) model, \Odrzygozdz\ showed that Property (T) does not hold at densities $d<\frac{3}{8}$ \cite{Odr16,Odr19}.
We can find a larger range of densities with Property FA but not Property (T) when considering the ``positive $\ell$-gonal model'', see Definition~\ref{def-l-gonal-positive-density}.
\begin{corollary}
	\label{cor-FA-no-T-range}
	For $\ell\in \{4,6,7,8,\ldots\}$ there exists $d_\ell > \frac{1}{\ell}$ so that
	for densities $\frac{1}{\ell}<d<d_\ell$ a.a.s.\ it holds that a random $\ell$-gonal group has Property FA but not Property (T).
	For $\ell \geq 5$, and densities $\frac{1}{\ell}<d<\frac{1}{4}$, a.a.s.\ it holds that a random positive $\ell$-gonal group has Property FA but not Property (T).
\end{corollary}
\begin{proof}
	Property FA follows from Theorems~\ref{thm-propFA-sharp}(3) and \ref{thm-positive-fa}.
	Known results about failure of Property (T) let us choose $d_\ell$ as follows.  For $\ell=4$, we may take $d_4 = \frac{3}{8}>\frac{1}{4}$ by \cite{Odr19}.
	For $\ell \geq 6$, we may take $d_\ell = \lfloor \ell/2 \rfloor / 2\ell \geq \frac{1}{5}>\frac{1}{\ell}$ by \cite[Theorem D]{Ashcroft22-14propT}.  

	We now show the claim for the positive $\ell$-gonal model.  At any density groups in the Gromov density model (at lengths which are multiples of $\ell$) have a finite index subgroup which is a quotient of a random group in the positive $\ell$-gonal model at the same density, see \cite[Lemma 3.15]{Kotowski13}.  Property (T) is preserved by quotients and finite over-groups, so the claim follows from the fact that in the Gromov density model a.a.s.\ property (T) fails at densities $d<1/4$ (\cite{Ashcroft22-14propT}, see also \cite{OllWise11, MacPrz15}).
\end{proof}
	
As already mentioned, the proof of Theorem~\ref{thm-propFA-sharp}(1,2) largely follows the arguments of \cite[Theorems 1 and 2]{AntLucSwi-13threshold}, who proved it in the case $\ell=3$, with some adaptations to the case of general $\ell$.
The proof of Theorem~\ref{thm-propFA-sharp}(3) is inspired by Dahmani--Guirardel--Przytycki's work in the Gromov density model \cite{DahGuiPrz11-nosplit}, but the specific arguments are different.  Their methods were adapted to that model, where $m$ is fixed and $\ell \to \infty$; the direct application of their methods (using their notion of `b-automata' etc.) gives Property FA for $d > 3/\ell$, whereas we can get the sharp density $d>1/\ell$.
In brief, one shows that for any potential splitting of $G = A \ast_C B$, the generators of $G$ have to have canonical forms that imply some relation has a non-trivial canonical form in $A \ast_C B$, giving a contradiction.

\subsection*{Outline of paper}
Theorem~\ref{thm-propFA-sharp} parts (1,2,3) are proved in Sections~\ref{sec-free}, \ref{sec-not-free-not-fa} and \ref{sec-fa} respectively.
Corollaries~\ref{cor-boundary-sharp} and \ref{cor-finite-out-sharp} are proved in Section~\ref{sec-boundaries}.

\subsection*{Acknowledgements}
The second named author thanks Mark Hagen for helpful comments.
We also thank the referees for their detailed and useful comments.

\section{Freeness}\label{sec-free}
	We consider $G \in \mathcal{M}(m,\ell,p)$ for $\ell \geq 3$ and $p \leq cm^{1-\ell}$ with $c>0$ small.
	First, we show there are a.a.s.\ at least two more generators than relators.
	\begin{lemma}
		\label{lem:RleqS-2}
		There exists $c>0$ such that for $p \leq cm^{1-\ell}$, $\langle S|R\rangle \in \mathcal{M}(m,\ell,p)$ satisfies $|R|\leq|S|-2$ a.a.s.
	\end{lemma}
	\begin{proof}
		There are \[ N=(2m-1)^\ell+1+(m-1)(1+(-1)^\ell) = 2^\ell(1+o(1))m^\ell \] cyclically reduced words of length $\ell$ in $m$ generators (see e.g.~\cite{Rivin10-growth-free-groups}).
		Let $X$ be the random variable counting the number of relations chosen; this has a binomial distribution $B(N,p)$ so 
		\[
			\mathbb{E}X = Np \leq c 2^\ell(1+o(1))m, \Var X = Np(1-p) = (1-o(1)) \mathbb{E}X.
		\]
		So by Chebyshev's inequality, assuming $c < \frac{1}{3} 2^{-\ell}$,
		\[
			P(|R| > |S|-2) = P(X > m-2) 
			\leq P( X > \mathbb{E}X + m^{1/2}(\Var X)^{1/2}) \leq \frac{1}{m} \to 0.
		\]
	\end{proof}

	We now state our key proposition.
	\begin{proposition} 
		[{cf.\ \cite[Lemma 4]{AntLucSwi-13threshold}}]
		\label{lem-delete-relator-generator}
		There exists $c>0$ such that for $p \leq cm^{1-\ell}$ the random group $\langle S | R \rangle \in \mathcal{M}(m,\ell,p)$ a.a.s.\ has the following property: for every non-empty subset $R' \subset R$ of relations there exists $a \in S, r \in R'$ such that neither $a$ nor $a^{-1}$ appears in any relation $t \in R' \setminus \{r\}$ and precisely one letter in $r$ belongs to the set $\{a,a^{-1}\}$.
	\end{proposition}
Before we prove Proposition~\ref{lem-delete-relator-generator}, we show how it applies to show freeness.
\begin{proof}
	[Proof of Theorem~\ref{thm-propFA-sharp}(1)]
	We fix $c>0$ small enough to satisfy both Lemma~\ref{lem:RleqS-2} and Proposition~\ref{lem-delete-relator-generator}, and now a.a.s.\ may assume that both their conclusions hold.
	Apply Proposition~\ref{lem-delete-relator-generator} with $R'=R$ to find $a \in S, r \in R$ so that $a$ or $a^{-1}$ appears exactly once in $r$, and not in any $t \in R\setminus\{r\}$.  Therefore one can delete $a$ from $S$ and $r$ from $R$ to see that $G = \langle S|R\rangle \cong \langle S\setminus\{a\} | R \setminus \{r\}\rangle$.  Now continue, applying Proposition~\ref{lem-delete-relator-generator} to $R' = R \setminus \{r\}$.  This process removes generators and relations together one by one.  By Lemma~\ref{lem:RleqS-2} $|S|-|R| \geq 2$ so this process ends showing $G$ is a free group with $|S|-|R|$ generators. 
\end{proof}
It remains to prove Proposition~\ref{lem-delete-relator-generator}, and this is the goal of the rest of this section.  In the proof we will use the following bound on component size in sparse random hypergraphs.
	\begin{theorem}
		[{Schmidt-Pruzan--Shamir~\cite[Theorem 3.6]{SchmidtPruzanShamir85}}]
		\label{thm-SPS-hypergraph-cpts}
		Let $\mathcal{H}$ be a random hypergraph with $m$ vertices and in which each subset with $\ell$ vertices appears independently with probability $\rho$.
		Suppose $(\ell-1)\rho\binom{m-1}{\ell-1} \leq \alpha$ for a constant $\alpha<1$.
		Then there exists a constant $K_\alpha$ such that with probability $1-o(m^{-2})$ the largest connected component in $\mathcal{H}$ contains at most $K_\alpha \log(m)$ vertices.
	\end{theorem}
% We now prove the proposition. %% maybe don't need? 
\begin{proof}
	[Proof of Proposition~\ref{lem-delete-relator-generator}]
	We follow the proof of \cite[Lemma 4]{AntLucSwi-13threshold}.

	For each relator $r \in R$, consider the number $\ell_r$ of different $s \in S$ 
	which appear in $r$ as $s$ or $s^{-1}$.   
	\begin{itemize}
		\item $r$ is of \emph{type 1} if $\ell_r \leq \ell-2$,
		\item $r$ is of \emph{type 2} if $\ell_r = \ell-1$, and
		\item $r$ is of \emph{type 3} if $\ell_r = \ell$.
	\end{itemize}
	For example, if $\ell=4$, $aba^{-1}b$, $aabc$, and $abcd$ are of types 1,2, and 3 respectively.
	Let $\mathcal{H}$ be the random hypergraph with vertex set $S$ and an edge for each $r \in R$ consisting of the generators $s \in S$ which appear in $r$ as $s$ or $s^{-1}$.
	Observe that a.a.s.\ each edge in $\mathcal{H}$ corresponds to a single relation $r \in R$.  This is because any subset of $S$ of size $\leq \ell$ could come from at most $O(1)$ different potential relations of length $\ell$, each of which is included in $R$ with probability $p$.  So the probability a given edge in $\mathcal{H}$ comes from two different relations in $R$ is $O(p^2)$, and summing over all possible edges we have that the probability $P$ that there exists an edge coming from at least two different relations satisfies:
	\begin{equation}\label{eq-no-double-edges}
		P \leq C_\ell m^\ell p^2 \leq C_\ell c^2 m^{2-\ell} \to 0,
	\end{equation}
	for some constant $C_\ell$ depending on $\ell$.  We use the convention that constants $C_\ell$ depend only on $\ell$, but may change from line to line.

	There are $\leq C_\ell m^{\ell-2}$ relators of type 1, for some constant $C_\ell$.  Let $X$ be the number of such relations.  Then
	\begin{equation}\label{eq-no-type-1}
	P(X>0) \leq \mathbb{E}X = C_\ell m^{\ell-2} p = cC_\ell m^{-1} \to 0, 
	\end{equation}
	so a.a.s.\ no relations of type 1 appear, so no edges of size $\leq \ell-2$ appear in $\mathcal{H}$.

	% We use the following result.
	Let $\mathcal{H}_\ell$ be the $\ell$-uniform hypergraph on vertex set $S$ consisting of all sets of $\ell$ generators which arise from a type 3 relation in $R$.
	The probability of each edge appearing is $\rho = 2^\ell \ell! p \leq 2^\ell \ell! c m^{1-\ell}$ so provided $c<(2^\ell \ell!(\ell-1))^{-1}$, Theorem~\ref{thm-SPS-hypergraph-cpts} applies to give that a.a.s.\ each connected component of $\mathcal{H}_\ell$ has $\leq K \log(m)$ vertices, where $K$ is a constant.

	We show each connected sub-hypergraph of $\mathcal{H}_\ell$ other than an isolated vertex, contains at least $2$ vertices belonging to exactly one edge of $\mathcal{H}_\ell$.
	If not, let $X$ be the number of non-trivial connected sub-hypergraphs on $k$ vertices, $\ell+1 \leq k \leq K\log(m)$, contained in $\mathcal{H}_\ell$ in which all but at most one vertex belongs to at least two edges.  
	Each such sub-hypergraph has at least
	$\lceil (2k-1)/\ell \rceil$ edges because each edge contains $\ell$ vertices.
	We use the following bound on this quantity.
	\begin{lemma}\label{lem:2k16k5} For $\ell \geq 3$ and $k \geq \ell+1$ we have $\lceil (2k-1)/\ell \rceil \geq 6k/5(\ell-1)$.
	\end{lemma}
	\begin{proof}
		Let $k=q\ell+r$ for some $q,r \in \mathbb{Z}$, $q \geq 1$, $0 \leq r \leq \ell-1$.

		If $r=0$ then $\lceil (2k-1)/\ell \rceil = 2q$, and $6k/5(\ell-1)=6q\ell/5(\ell-1) \leq 18q/10 \leq 2q$.
		If $1\leq r \leq \frac{1}{2}(\ell+1)$ then 
		\begin{align*}
			\frac{k/(\ell-1)}{\lceil (2k-1)/\ell \rceil} &
			= \frac{q+(q+r)/(\ell-1)}{2q+1}
			\leq \frac{q+(q+\frac{1}{2}(\ell+1))/(\ell-1)}{2q+1}
			\\ & \leq \frac{1+(1+\frac{1}{2}(3+1))/(3-1)}{3}
			= \frac{5}{6},
		\end{align*}
		since the worst cases are when $q \geq 1$ and then $\ell \geq 3$ are minimized.
		Similarly, if $\frac{1}{2}(\ell+1) < r \leq \ell-1$ then
		\begin{align*}
			\frac{k/(\ell-1)}{\lceil (2k-1)/\ell \rceil}
			& = \frac{q+(q+r)/(\ell-1)}{2q+2}
			\leq \frac{q+(q+(\ell-1))/(\ell-1)}{2q+2}
			\\ & \leq \frac{1+(1+(3-1))/(3-1)}{4}
			= \frac{5}{8} \leq \frac{5}{6}. \qedhere
		\end{align*}
	\end{proof}

	Let $Y$ be the number of sub-hypergraphs of $\mathcal{H}_\ell$ on $k$ vertices, $\ell+1\leq k \leq K\log(m)$, with exactly $\lceil 6k/5(\ell-1) \rceil$ edges.
	For a given set of $k$ vertices, the probability there is a subgraph of $\mathcal{H}_\ell$ consisting of these $k$ vertices and having exactly $\lceil 6k/5(\ell-1) \rceil$ edges is at most
	\begin{align*}
		\binom{\binom{k}{\ell}}{\lceil 6k/5(\ell-1) \rceil} \rho^{\lceil 6k/5(\ell-1) \rceil}
		& \leq ( k^\ell)^{\lceil 6k/5(\ell-1) \rceil} \rho^{\lceil 6k/5(\ell-1) \rceil}
		\\ & \leq (2^\ell \ell!  k^\ell p)^{6k/5(\ell-1)},
	\end{align*}
	so, letting $C_\ell$ be a constant which may change from line to line, 
	\begin{align*}
		\mathbb{E}Y & \leq \sum_{k=\ell+1}^{K\log m} \binom{m}{k} (2^\ell\ell!k^\ell p)^{6k/5(\ell-1)}
		\leq \sum_{k=\ell+1}^{K\log m} \left(\frac{em}{k}\right)^k \left(\frac{2^\ell\ell!k^\ell c}{m^{\ell-1}}\right)^{6k/5(\ell-1)}
		\\ & \leq \sum_{k=\ell+1}^{K\log m} \left(\frac{C_\ell k^{(\ell+5)/(\ell-1)}}{m}\right)^{k/5}
		\leq K \log(m) \left( \frac{C_\ell (K \log(m))^{4}}{m}\right)^{4/5} \to 0.
	\end{align*}
	Since 
	\begin{equation}\label{eq-each-component-two-exposed}
		P(X>0) \leq P(Y>0) \leq \mathbb{E}Y,
	\end{equation}
	a.a.s.\ it holds that each connected sub-hypergraph of $\mathcal{H}_\ell$ other than an isolated vertex, contains at least $\ell-1$ vertices belonging to exactly one edge of $\mathcal{H}_\ell$.

	Let us now consider how type 2 edges behave.
	First, a.a.s.\ each type 2 edge in $\mathcal{H}$ meets each component of $\mathcal{H}_\ell$ in at most one vertex.
	If $Z$ counts the number of components of $\mathcal{H}_\ell$ which meet a type 2 edge in at least two vertices, since we may assume each connected component of $\mathcal{H}_\ell$ has at most $K \log(m)$ vertices, there are at most $C_\ell (K \log m)^2 m^{\ell-3}$ choices of potential type 2 edges to consider meeting each connected component.
	So we have,
	\begin{equation}\label{eq-component-meets-type-2-at-1-vertex}
		P(Z>0) \leq \mathbb{E}Z \leq C_\ell m (K \log m)^2 m^{\ell-3} p = C_\ell (\log m)^2 m^{1+\ell-3+1-\ell} \to 0.
	\end{equation}
	Second, a.a.s.\ no component of $\mathcal{H}_\ell$ meets two type 2 edges.
	For if $Z'$ counts the number of components of $\mathcal{H}_\ell$ which meet at least two type 2 edges, for a given component there are at most $C_\ell (K \log m)^2$ choices of two vertices contained in a type 2 edge, and at most $C_\ell m^{2(\ell-2)}$ choices of the remaining vertices in each edge, therefore,
	\begin{align}
		P(Z'>0) & \leq \mathbb{E}Z' \leq C_\ell m (K \log m)^2 m^{2(\ell-2)}p^2 
		\notag \\ & =C_\ell (\log m)^2 m^{1+2\ell-4+2-2\ell} \to 0.
		\label{eq-component-meets-at-most-one-type-2}
	\end{align}
	Third, a.a.s.\ no two type 2 edges share a vertex, i.e.\ the type 2 edges form a matching.
	For if $Z''$ counts the number of vertices in $S$ which lie in at least two type 2 edges,
	\begin{equation}\label{eq-type-2-matching}
		P(Z''>0) \leq \mathbb{E}Z'' \leq C_\ell m m^{2(\ell-2)}p^2 = C_\ell m^{1+2\ell-4+2-2\ell} \to 0.
	\end{equation}

	Lemma~\ref{lem-delete-relator-generator} now follows:
	a.a.s.\ we may assume that the conditions of \eqref{eq-no-double-edges}, \eqref{eq-no-type-1}, \eqref{eq-each-component-two-exposed}, \eqref{eq-component-meets-type-2-at-1-vertex}, \eqref{eq-component-meets-at-most-one-type-2} all hold. 
	Take any subset $R' \subset R$ of relations, and consider the corresponding sub-hypergraph $F$ of $\mathcal{H}$.
	By \eqref{eq-no-double-edges}, each edge in $\mathcal{H}$ corresponds to a unique relation in $R$, and by \eqref{eq-no-type-1} there are no type 1 edges.
	Let $F' = F \cap \mathcal{H}_\ell$.
	If $F'$ has an edge then by \eqref{eq-each-component-two-exposed} each of its connected components contains at least two vertices not in any other edge of $F'$, and by \eqref{eq-component-meets-type-2-at-1-vertex} and \eqref{eq-component-meets-at-most-one-type-2} at least one of these vertices is in no other hyperedge, so this vertex $a$ and the relation $r$ it lies in satisfy the conclusion of the theorem.

	If $F'$ has no edges, then $F$ consists entirely of type 2 edges, which by \eqref{eq-type-2-matching} form a matching.  So choose $r \in R'$ corresponding to one such edge.  As it is type 2, there are $\ell-2 \geq 1$ choices of $a$ in $r$ which satisfy the conclusion of the lemma.
\end{proof}

\section{Not free and not FA}\label{sec-not-free-not-fa}

Ashcroft and Roney-Dougal show hyperbolicity of random groups in the $\mathcal{M}(m,\ell,d)$ model at densities $d<1/2$, where they choose as relations a set of cyclically reduced words of length $\ell$ of size $\lfloor (2m-1)^{\ell d} \rfloor$ uniformly at random from amongst all such sets.  They show, building on work of Ollivier (see \cite[Chapter V]{Oll05} and references therein):
\begin{theorem}
	[{\cite[Theorem 3.11]{AshcroftRonDou20}}]
	\label{thm-uniform-hyperbolic}
	For any $\ell \geq 3, d_* < 1/2, \epsilon>0$, there exists $f:\N\to [0,1]$ with $f(m) \to 0$ as $m\to\infty$ so that for any $0< d < d_*$ with probability $\geq 1-f(m)$  a random group presentation $G \in \mathcal{M}(m,\ell,d)$ 
	has the property that any reduced van Kampen diagram $D$ for $G$ satisfies the linear isoperimetric inequality $|\partial D| \geq \ell(1-2d-\epsilon) |D|$.

	In particular, $G=\langle S|R \rangle$ is a.a.s.\ infinite torsion-free hyperbolic with aspherical presentation and Euler characteristic $\chi(G)=1-|S|+|R|$.
\end{theorem}
This theorem is stated slightly differently in \cite{AshcroftRonDou20}, but their proof gives the uniform probability bounds listed above.  We refer to their paper for definitions which we do not need in this paper.
The following standard argument connects the two models $\mathcal{M}(m,\ell,d)$ and $\mathcal{M}(m,\ell,p)$ for $p=m^{\ell(1-d)}$.
\begin{proposition}\label{prop-uniform-to-binomial}
	Let $\mathcal{P}$ be a property of a group presentation.
	Fix $\ell \geq 3$.
	Suppose there exists $d_0,d_1 \in [0,1]$ and $f(m) \to 0$ as $m \to \infty$ so that for $d \in (d_0,d_1)$ the probability $G \in \mathcal{M}(m,\ell,d)$ has $\mathcal{P}$ is at least $1-f(m)$.
	Then for any $\delta>0$, and $p=p(m) \in [m^{\ell(d_0+\delta-1)},m^{\ell(d_1-\delta-1)}]$, a.a.s.\ it holds that $G \in \mathcal{M}(m,\ell,p)$ satisfies $\mathcal{P}$.
\end{proposition}
\begin{proof}
	This follows from \cite[Proposition 1.12]{JansonLuczakRucinski-00-random-graphs-book}: this is applied with $N=2^\ell (1+o(1))m^\ell$, $p(m) \in [m^{\ell(d_0+\delta-1)},m^{\ell(d_1-\delta-1)}]$, $q(m)=1-p(m)=1-o(1)$.
	Then, $ 2^\ell(1-o(1))m^{\ell (d_0 + \delta)} \leq Np \leq 2^\ell(1+o(1))m^{\ell (d_1 - \delta)} $, so $Np \to \infty$ as $m \to \infty$.
	Therefore for any sequence $M(m)=Np+O(\sqrt{Npq})$,
	\begin{align*}
		M(m) = Np + O(\sqrt{Npq}) &= Np \left( 1 + O \left( \frac{1}{\sqrt{Npq}} \right) \right) \\ 
		& \in \left[ \frac{Np}{2}, 2Np \right] 
						   \subset \left[ m^{\ell ( d_0 + \frac{\delta}{2} )}, m^{\ell ( d_1 - \frac{\delta}{2} )} \right] 
						   \\ & \subset \left( (2m-1)^{\ell  d_0 }, (2m-1)^{\ell d_1 } \right) . \qedhere
	\end{align*}
\end{proof}
There is a converse statement for monotone properties, see \cite[Proposition 1.13]{JansonLuczakRucinski-00-random-graphs-book}, but we do not require that in this paper. 
\begin{corollary}
	\label{cor-binomial-hyperbolicity}
	For any $\ell \geq 3, \epsilon>0$, 
	if $p< m^{-\ell/2-\epsilon}$ then
	a.a.s.\ a random group presentation $G \in \mathcal{M}(m,\ell,p)$ 
	has the property that any reduced van Kampen diagram $D$ for $G$ satisfies the linear isoperimetric inequality $|\partial D| \geq \frac{1}{2}\epsilon |D|$.

	In particular, $G=\langle S|R \rangle$ is a.a.s.\ infinite torsion-free hyperbolic with aspherical presentation and Euler characteristic $\chi(G)=1-|S|+|R|$.
\end{corollary}
\begin{proof}
	This follows from Proposition~\ref{prop-uniform-to-binomial} with $d_0=0, d_1=\frac{1}{2}-\frac{\epsilon}{2\ell}, \delta=\frac{\epsilon}{2\ell}$, and Theorem~\ref{thm-uniform-hyperbolic} with $d_*=d_1$, `$\epsilon$'$=\frac{\epsilon}{2\ell}$.
\end{proof}
This corollary could alternatively have been proved directly by running through the arguments of~\cite{AshcroftRonDou20}.

Now we turn to the goal of this section: to prove Theorem~\ref{thm-propFA-sharp}(2).
First, let us prove that there are at least $3m$ relations in $R$ for $p$ large enough.
	\begin{lemma}[{Compare \cite[Lemma 11]{AntLucSwi-13threshold}}]
		\label{lem-lots-of-relations}
		If $p \geq 4 m^{1-\ell}$ then a.a.s.\ $G=\langle S|R \rangle \in \mathcal{M}(m,\ell,p)$ has $|R| \geq 3m$.
	\end{lemma}
	\begin{proof}
	Let us denote by $X$ the random variable counting the number of relations containing $\ell$ different elements in $G$.
	These relations are picked among a number $N$ of relations, 
	with $N = 2m \cdot (2m-2) \cdots (2m-2(\ell-1)) = 
	\frac{m!}{(m-\ell)!} 2^{\ell} $.
		Let us compute the variance and expectation of $X$, assuming as we may that $m \geq 2\ell$ and $p=o(1)$:
	\begin{align*}
		\Var X & = N p ( 1 - p ) = ( 1-p ) EX = ( 1 + o(1) ) \mathbb{E}X, \text{ and} \\
		\mathbb{E}X & = Np = \frac{m!}{(m-\ell)!} 2^{\ell} p 
		\geq \frac{m\cdots(m-\ell+1)}{m^\ell} 2^\ell \cdot 4m \geq 4m.
	\end{align*}
		By Chebyshev, 
		\[
			P(|R| < 3m) \leq P(X < \mathbb{E}X - (1/4)\mathbb{E}X)
			\leq \frac{16\Var X}{(\mathbb{E}X)^2} \leq \frac{5}{m} \to 0. \qedhere
		\]
	\end{proof}
	
Now let us show that for $p$ not too large, a.a.s.\ some generators do not appear in any relator.
	\begin{lemma}[{Compare \cite[Lemma 12]{AntLucSwi-13threshold}}]
		\label{lem-missing-generators}
		If $p \leq \ell^{-1}2^{-\ell-2} \log(m) m^{1-\ell}$ then a.a.s.\ at least $m^{1/2}/2$ of the generators of $G = \langle S|R \rangle \in \mathcal{M}(m,\ell,p)$ do not appear in any $r \in R$, as themselves or their inverses.
	\end{lemma}
	\begin{proof}
		We roughly follow the strategy of \cite{AntLucSwi-13threshold}, but add many details.
		For $s \in S$, let $Y_s=1$ if neither $s$ nor $s^{-1}$ appears in any $r \in R$, and let $Y_s=0$ otherwise.  Let $Y = \sum_{s \in S} Y_s$ count the number of generators which do not appear in any relation.

		Let $N_s$ count the number of cyclically reduced words of length $\ell$ in $S$ which contain $s$ or $s^{-1}$.  
		We can bound $N_s$ from above by choosing a location for an $s$ or $s^{-1}$, then filling in the remaining $\ell-1$ letters each with $2m$ options.  
		On the other hand, $N_s$ is bounded below by choosing $s$ or $s^{-1}$ for the first letter, there are $(2m-1)$ choices for the next letters until the $\ell$th letter has at least $(2m-2)$ choices.  So if $m$ is large enough depending only on $\ell$ we have:
		\begin{equation}
			\label{eq:Ns}
			2^{\ell-1}m^{\ell-1} \leq 2 (2m-1)^{\ell-2}(2m-2) 
			\leq N_s \leq 2\ell(2m)^{\ell-1}=2^\ell\ell m^{\ell-1}.			
		\end{equation}
		Note that $N_s$ does not depend on $s$, and that 
		\[
			q = P(Y_s=1) = (1-p)^{N_s} 
			\geq e^{-2p N_s}
			\geq e^{-\log(m)/2} = \frac{1}{\sqrt{m}}.
		\]
		Therefore,
		\[
			\mathbb{E}Y = mq \geq \sqrt{m} \to \infty \text{ as } m \to \infty.
		\]

		Now for $s \neq t \in S$, let $Z_{s,t}$ be the indicator for the event that neither $s$ nor $t$ (nor their inverses) appear in any relator in $R$.
		Let $N_{s,t}$ be the number of cyclically reduced words of length $\ell$ in $S$ which contain both $s$ and $t$.
		Observe that, as we can choose two locations for an $s$ or $s^{-1}$ and for a $t$ or $t^{-1}$, then the remaining $\ell-2$ locations freely, we have
		\[
			N_{s,t} \leq (2\ell)^2(2m)^{\ell-2} \leq \frac{2\ell^2}{m}N_s,
		\]
		using \eqref{eq:Ns}.
		Thus as $Y_sY_t$ indicates the event that none of the $N_s+N_t-N_{s,t}$ relations containing either $s$ or $t$ are in $R$, we have the covariance of $Y_s,Y_t$ satisfies:
		\begin{align*}
			\mathrm{Cov}(Y_s,Y_t) & = \mathbb{E}(Y_sY_t) - \mathbb{E}(Y_s)\mathbb{E}(Y_t)
			\\ & = (1-p)^{N_s+N_t-N_{s,t}} - q^2
			\\ & = q^2 \left( \frac{1}{(1-p)^{N_{s,t}}} -1 \right)
				\leq q^2 \left( \frac{1}{q^{2\ell^2/m}}-1\right)
				\leq q^2 \left( m^{\ell^2/m} -1 \right).
		\end{align*}
		Therefore,
		\begin{align*}
			\Var(Y) & = \sum_{s \in S} \Var(Y_s) + \sum_{s \neq t \in S} \mathrm{Cov}(Y_s,Y_t)
			\\ & \leq mq(1-q) + m^2 q^2 \left( m^{\ell^2/m}-1\right)
			\\ & \leq \mathbb{E}Y + (\mathbb{E}Y)^2 \left( m^{\ell^2/m} -1 \right)
			= o((\mathbb{E}Y)^2)
		\end{align*}
		as $\mathbb{E}Y = mq \to \infty$.
		So by Chebyshev, 
		\[
			P\Big(Y \leq \frac{1}{2}\sqrt{m}\Big) 
			\leq P\Big(Y \leq \mathbb{E}Y-\frac{1}{2}\mathbb{E}Y\Big)
			\leq \frac{4\Var(Y)}{(\mathbb{E}Y)^2} = o(1).\qedhere
		\]
	\end{proof}
We now can find a range of $p$ for which a.a.s.\ the group is neither free nor has property FA.
\begin{proof}
	[Proof of Theorem~\ref{thm-propFA-sharp}(2)]
	Let $p$ satisfy $4 m^{1-\ell} \leq p \leq \ell^{-1}2^{-\ell-2} \log(m)m^{1-\ell}$.
	We may assume a.a.s.\ that $G=\langle S|R \rangle \in \mathcal{M}(m,\ell,p)$ satisfies the conclusions of Corollary~\ref{cor-binomial-hyperbolicity} and Lemmas~\ref{lem-lots-of-relations} and \ref{lem-missing-generators}.

	By Corollary~\ref{cor-binomial-hyperbolicity} we have that the presentation complex for $G=\langle S | R \rangle$ is aspherical.
	Therefore the Euler characteristic $\chi(G) = 1-|S|+|R| \geq 1-m+3m>0$ by Lemma~\ref{lem-lots-of-relations}, so $G$ is not free.

	By Lemma~\ref{lem-missing-generators}, at least $\frac{1}{2}\sqrt{m}$ generators do not appear in any relation of $G$, so $G \cong H \ast F_q$ for some $q \geq \frac{1}{2}\sqrt{m}$ where $H$ is the subgroup of $G$ generated by the other generators.
\end{proof}

\section{Property FA}\label{sec-fa}
We can show FA using only relators which contain positive words.
\begin{definition}
	\label{def-l-gonal-positive-density}
	Let $d \in (0,1)$ be a density, and $\ell \in \{3,4,\ldots\}$ a fixed length.
	A \emph{(binomial) random positive $\ell$-gonal group with probability $p$ (and $m \geq 2$ generators)} $\mathcal{P}(m,\ell,p)$ is the probability space consisting of group presentations $G = \langle S | R \rangle$ where $S = \{s_1,\ldots, s_m\}$, and $R$ consists of relators where each of the $m^\ell$ words of length $\ell$ in $S$ (not $S^{-1}$) are chosen independently with probability $p$.  
\end{definition}
Observe that groups in $\mathcal{M}(m,\ell,p)$ are quotients of groups in $\mathcal{P}(m,\ell,p)$ for the same values of $m,\ell,p$.  Since Property FA is preserved by quotients, to prove Theorem~\ref{thm-propFA-sharp}(3) it suffices to show:
\begin{theorem}
	\label{thm-positive-fa}
	Fix $\ell \geq 3$.  There exists $C'>0$ so that for $p \geq C' \log(m)m^{1-\ell}$,
	a.a.s.\ $G \in \mathcal{P}(m,\ell,p)$ satisfies Property FA.
\end{theorem}
Before proving the theorem, we establish some preliminary lemmas about the a.a.s.\ guaranteed existence of certain relators.
\begin{lemma}
	\label{lem-exist-VVV-words}
	Fix $\ell \geq 3$.  For any $C'>0$ so that $p \geq C' \log(m)m^{1-\ell}$,
	a.a.s.\ $G = \langle S|R \rangle \in \mathcal{P}(m,\ell,p)$ satisfies:

	\begin{enumerate}
		\item[(L)]
			\makeatletter\def\@currentlabel{(L)}\makeatother
			\label{propL}
	For any subsets $V_1, V_2, \ldots, V_\ell \subset S$ with each $|V_i| \geq m/100$, there exists a relator $r \in R$ of the form $r=v_1v_2\cdots v_\ell$ with each $v_i\in V_i$.
	\end{enumerate}
\end{lemma}
Here `\ref{propL}' refers to relators with letters coming from `Large' subsets of $S$.
\begin{proof}
	Let us denote $P_L$ the probability that `\ref{propL}' is satisfied.
	For each $i$, there are $\binom{m}{|V_i|} \leq 2^m$ choices for $V_i$, so there are at most $(2^m)^{\ell}$ sequences $V_1, \cdots, V_{\ell}$ satisfying $\left| V_i \right| \geq m / 100 $ for all $i$. If the claim fails for one of these sequences, each of the $\geq (m/100)^\ell$ words of the prescribed form is omitted from $R$. Using the bound $1+x \leq e^x$, we find that the claim fails with probability:
	\begin{align*}
		1-P_{L} & \leq (2^m)^\ell (1-p)^{(m/100)^\ell}
		\leq 2^{m\ell} \exp\left(-C'\log(m)m^{1-\ell}100^{-\ell}m^\ell\right)
		\\ & \leq \exp\left(m\ell\log(2) -C'100^{-\ell} \log(m)m\right) 
		\to 0
	\end{align*}
	as $m \to \infty$.
\end{proof}
Slightly more technical is the following.
\begin{lemma}
	\label{lem-exist-UVV-words}
	Fix $\ell \geq 3$.  There exists $C'>0$ so that for $p \geq C' \log(m)m^{1-\ell}$,
	a.a.s.\ $G = \langle S|R \rangle \in \mathcal{P}(m,\ell,p)$ satisfies:
	\begin{enumerate}
		\item[(SL)]
			\makeatletter\def\@currentlabel{(SL)}\makeatother
			\label{propSL}
	For any partition $S = U \sqcup V$ with $1 \leq |U| \leq 99m/100$, there exists a relator $r \in R$ of the form $r=uv_1v_2\cdots v_{\ell-1}$ with $u \in U$ and $v_i \in V$ for $1 \leq i \leq \ell-1$.
	\end{enumerate}
\end{lemma}
Here `\ref{propSL}' refers to relators with initial letter coming from a potentially `Small' subset of $S$, and all remaining letters coming from a `Large' subset of $S$.
\begin{proof}
	We first bound the probability the claim fails for some choice of $U$, with $k=|U|$ satisfying $1 \leq k \leq 99m/100$.
	There are $\geq k (m/100)^{\ell-1}$ words consisting of a letter from $U$ followed by $\ell-1$ letters from $V$.
	The probability, denoted $P_R$, that $R$ misses all of these is
		\begin{align*}
			P_R & \leq (1-p)^{k (m/100)^{\ell-1}} 
			 \leq \exp\left( -C'\log(m)m^{1-\ell} k 100^{1-\ell}m^{\ell-1} \right)
			\\  & \leq \exp\left( -C'100^{-\ell} k \log(m) \right).
		\end{align*}
	Summing over the choices of $k$ and $U$ (which determines $V$), the probability, denoted $P$, that the lemma fails is at most:
		\begin{equation*}
			P \leq \sum_{k=1}^{\lfloor 99m/100 \rfloor} \binom{m}{k} e^{-C' 100^{-\ell} k \log(m)}
			 \leq \sum_{k=1}^{\lfloor 99m/100 \rfloor} m^k m^{-C' 100^{-\ell} k}.
		%\label{eq-few-positive}
		\end{equation*}
		Provided $C'> 2\cdot 100^\ell$, this sum is
		\[
			\leq \sum_{k=1}^\infty m^{-k} \leq \frac{2}{m} \to 0
		\]
		as $m \to \infty$. 
\end{proof}
The properties \ref{propL} and \ref{propSL} are what we use to show $G$ has property FA.  First we show there is no surjection onto $\mathbb{Z}$.
	\begin{lemma}
		\label{lem-no-surject-onto-Z}
	Let $G = \langle S|R \rangle$ with $|S|=m$ and each $r \in R$ a (positive) word of length $\ell$ in $S$.
		Then if this presentation satisfies \ref{propL} and \ref{propSL}, $G$ does not surject onto $\mathbb{Z}$.
	\end{lemma}
	\begin{proof}
		Suppose there exists a surjection $\phi:G \to \mathbb{Z}$.
		Let $S = P \sqcup N \sqcup Z$ be the partition of $S$ according to whether $\phi(s)$ is positive, negative or zero.  
		Up to replacing $\phi$ by $-\phi$, we may assume $|P| \geq |N|$.
		If $|P| \geq m/100$ then \ref{propL} applied with all $V_i=P$ gives that there exists a relator $r$ consisting of letters from $P$, but then $\phi(r)>0$, contradiction.
		So $|P| < m/100$, and hence $|N| < m/100$, thus $|Z| > 98m/100$.
		We must have $P$ and possibly $N$ are non-empty, else $\phi$ is not surjective.
		Then by \ref{propSL} applied with $U = P \cup N$ and $V = Z$ there exists a relator $r=uv_1\cdots v_{\ell-1}$ with $u \in P \cup N$ and each $v_i \in Z$, so $\phi(r)=\phi(u) \neq 0$, contradiction.
		So no such $\phi$ exists.
	\end{proof}
Now we establish property FA.
\begin{proposition}\label{prop-L-SL-has-FA}
	Let $G = \langle S|R \rangle$ with $|S|=m$ and each $r \in R$ a (positive) word of length $\ell$ in $S$.
		Then if this presentation satisfies \ref{propL} and \ref{propSL}, $G$ has Property FA.
\end{proposition}

\begin{proof}[Proof of Proposition~\ref{prop-L-SL-has-FA}]
	By \cite[Theorem I.15]{SerreTrees03}, the finitely generated group $G$ has property FA if and only if it does not surject onto $\mathbb{Z}$ and is not an amalgamated product $G = A \ast_C B$.  A surjection onto $\mathbb{Z}$ is ruled out by Lemma~\ref{lem-no-surject-onto-Z}.  
	%Throughout, we assume that the hypotheses of Proposition~\ref{prop-L-SL-has-FA} hold.

	Now we consider whether $G$ can be written as an amalgamated product $G = A \ast_C B$.
	The overall strategy is similar in spirit to that of \cite{DahGuiPrz11-nosplit}.  We show that if $G$ can be written as an amalgamated product $G=A \ast_C B$ then by considering the possible behaviour of reduced forms for all the generators in $S$ one derives a contradiction by showing some relator is non-trivial. We use \ref{propL} and \ref{propSL}, whereas Dahmani--Guirardel--Przytycki use large automata.

	Let $G = A \ast_C B$.  Every $g \in G$ can be written in \emph{reduced form} where:
	\begin{itemize} 
		\item $g \in G\setminus C$ is written as $g=a_1b_1a_2b_2\cdots a_kb_k$ for some $a_i \in A, b_i \in B$, so that none of $b_1,a_2,\ldots,a_k$ are in $C$, and $a_1,b_k$ are not in $C\setminus \{e\}$. 
		\item $g \in C$ is written as $g=ge$ with $k=1, a_1=g, b_1=e$. 
	\end{itemize}
	The reduced form is unique up to choosing suitable transversals for $C$ in $G$, see Serre~\cite[Section I.1.2]{SerreTrees03}.  In particular, for $g \in G$ the left coset $[a_1] \in A/C$ is uniquely determined by $g$.
	The \emph{length} of $g$ with respect to the splitting $G = A \ast_C B$ is the number of non-trivial terms in the reduced form of $g$.
	The uniqueness of reduced forms implies that length is also well-defined; in particular if $g$ has a reduced form with non-zero length, then $g \neq e$ in $G$ 

	Suppose $G$ splits as an amalgamated product $G=A \ast_C B$.
	We can and do assume that the total length of elements of $S$ with respect to the splitting is minimal among all such ways to write $G$ as an amalgamated product.  This assumption will be needed after Lemma~\ref{lem:AB-BA-BB-small}.

	Let $\cA, \cB, \cC, \cD$ denote the partition of $S$ into subsets which consist of generators which lie in $A\setminus C, B\setminus C, C, G \setminus (A \cup B)$ respectively.
	Let $\alpha,\beta,\gamma,\delta$ denote the sizes of $\cA, \cB, \cC, \cD$ respectively.
	\begin{lemma}[{cf.\ \cite[Lemma 2.10]{DahGuiPrz11-nosplit}}] \label{lem:alpha-beta-gamma-small}
		For any non-trivial splitting of $G$ as an amalgamated product, $\alpha+\beta+\gamma \leq m/50$.
	\end{lemma}
	\begin{proof}
		Suppose $\alpha+\gamma \geq m/100$.
		Let $U = \cB \cup \cD$, $V = \cA \cup \cC$; note that $U \neq \emptyset$ else $G=A$.
		Then \ref{propSL} gives that there exists a relator
		$r = u v_1\cdots v_{\ell-1}$ with $u \in \cB \cup \cD$ so $u \notin A$ and each $v_i \in \cA \cup \cC \subset A$.  But then in $G$, $u = (v_1\cdots v_{\ell-1})^{-1} \in A$, contradiction.
		So $\alpha +\gamma < m/100$.

		Similarly, $\beta+\gamma < m/100$, else $G=B$.  So the lemma is proved.
	\end{proof}

	Elements of $\cD$ have reduced forms of lengths $\geq 2$.
	Let $\cD_{AA}, \cD_{AB}$, $\cD_{BA}$, $\cD_{BB}$ partition $\cD$ according to the first and last non-trivial letters of the reduced form, e.g., if $s \in \cD$ has reduced form starting with an element of $A\setminus C$ and ending with an element of $B\setminus C$, we have $s \in \cD_{AB}$.  Let $\delta_{AA} = |\cD_{AA}|$, etc.
	We can derive some immediate bounds on the sizes of these sets.
	\begin{lemma} \label{lem:AB-BA-BB-small}
		For any non-trivial splitting of $G$ as an amalgamated product, 
		$\delta_{AB} < m/100, \delta_{BA} < m/100$, and either $\delta_{AA}<m/100$ or $\delta_{BB} < m/100$.
	\end{lemma}
	\begin{proof}
		If $\delta_{AB} \geq m/100$, then by \ref{propL} with each $V_i = \cD_{AB}$, there exists $r \in R$ with $r=v_1v_2\cdots v_\ell$, with each $v_i \in \cD_{AB}$.  But then concatenating the reduced forms for the $v_i$ gives a reduced form for $r$ with respect to the splitting, so the length of $r$ is positive, contradiction.

		Likewise, $\delta_{BA} \geq m/100$ leads to a contradiction.

		If both $\delta_{AA} \geq m/100$ and $\delta_{BB} \geq m/100$, then \ref{propL} applied with $V_i = \cD_{AA}$ for $i$ odd and $V_i=\cD_{BB}$ for $i$ even leads to a contradiction.
	\end{proof}
	Up to swapping $A$ and $B$, we have reduced to the case that
	$\cD_{AA}$ has size at least $95m/100$.
	Now we have to further partition $\cD_{AA}$ by considering the first and last letters of reduced forms of letters in $\cD_{AA}$, and use our assumption that the splitting is chosen so that the total lengths of generators is minimal.
	
	For each left coset $[a] \in A/C$, let $N_{[a]}$ denote the number of reduced forms of letters in $\cD_{AA}^{\pm}$ which begin with an element of $[a]$, so $0 \leq N_{[a]} \leq 2m$.   This counts the number of $s \in \cD_{AA}$ which have a reduced form starting with an element of $aC$, plus the number of $s \in \cD_{AA}$ which have a reduced form ending with an element of $Ca^{-1}$ (as the reduced form of $s^{-1}$ will then start with an element of $aC$).  Minimality lets us bound $N_{[a]}$ in a much stronger way:
	\begin{lemma}
		For each $[a] \in A/C$, $N_{[a]} \leq m/10$.
	\end{lemma}
	\begin{proof}
		We consider what happens when we conjugate the splitting $G = A \ast_C B$ by $a$.
		For each $s \in S \setminus \cD_{AA}$, the reduced form of $s$ increases in length by at most two, a total increase in length of $\leq 10m/100$.
		Now consider $s \in \cD_{AA}$, with reduced form $s=a_1b_1\cdots a_k$.
		The reduced form of $a^{-1}sa$ is $(a^{-1}a_1)b_1 \cdots b_{k-1}(a_ka)$, up to potentially absorbing $(a^{-1}a_1)$ and/or $(a_ka)$ into the adjacent terms if they are in $C$.  But summing over each $s \in \cD_{AA}$, exactly $N_{[a]}$ of the start and end terms of the reduced forms lie in $C$.
		Thus the total change in reduced lengths for the generators in $\cD_{AA}$ is $-N_{[a]}$.
	
		Since the splitting $G=A\ast_C B$ was chosen to have minimal total length of reduced forms for generators, we have $N_{[a]} \leq 10m/100$.
	\end{proof}
	Since each $N_{[a]}$ is at most $m/10$, and 
	\[ \sum_{[a] \in A/C} N_{[a]} = 2|\cD_{AA}| \in [19m/10,2m],\]
	we can write $A/C = X \sqcup Y$ so that
	$\sum_{[a] \in X} N_{[a]}$ and $\sum_{[b] \in Y} N_{[b]}$ both lie between $9m/10$ and $11m/10$.
	We now run a similar argument for the last time, partitioning $\cD_{AA}$ as $\cD_{XX} \sqcup \cD_{XY}\sqcup \cD_{YX} \sqcup \cD_{YY}$, where e.g.\ $s \in \cD_{XY}$ indicates that the first letter in the reduced form for $s$ lies in some $[a] \in X$, and that the first letter in the reduced form for $s^{-1}$ lies in some $[a] \in Y$.

	\begin{lemma}
		Both $|\cD_{XY}| < m/100$ and $|\cD_{YX}|<m/100$.
	\end{lemma}
	\begin{proof}
		Suppose $|\cD_{XY}| \geq m/100$; the other case is similar.
		If $s,s' \in \cD_{XY}$ have reduced forms $s = a_1b_1 \cdots b_{k-1}a_k$ and $s' = a_1'b_1'\cdots a_{k'}'$, then a reduced form for $ss'$ is
		\[
			ss' = a_1b_1\cdots b_{k-1} (a_ka_1') b_1' \cdots a_{k'}'.
		\]
		This is because $a_ka_1' \in C$ if and only if $a_1'C = (a_k)^{-1}C$, but $a_1'C \in X$ and $(a_k)^{-1}C \in Y$.
		Thus letting every $V_i = \cD_{XY}$, by \ref{propL} we find a relator $r = v_1\cdots v_\ell$ with each $v_i \in \cD_{XY}$, which has positive length with respect to the splitting, contradiction.
	\end{proof}
	Now $2|\cD_{XX}| \leq \sum_{[a]\in X} N_{[a]} \leq 11m/10$ so $|\cD_{XX}| \leq 11m/20$.
	Likewise $|\cD_{YY}| \leq 11m/20$.
	But $|\cD_{XX}|+|\cD_{YY}| \geq 93m/100$, so both $|\cD_{XX}|$ and $|\cD_{YY}|$ are $ \geq 38m/100$.
	Apply \ref{propL} with $V_i = \cD_{XX}$ for $i$ odd and $V_i = \cD_{YY}$ for $i$ even to derive the final contradiction.
\end{proof}

We may now complete the proof of Theorem~\ref{thm-positive-fa}.
\begin{proof}
	[Proof of Theorem~\ref{thm-positive-fa}]
	Fix the choice of $C'$ from Lemma \ref{lem-exist-UVV-words}.
	The theorem now follows from Lemmas~\ref{lem-exist-VVV-words} and \ref{lem-exist-UVV-words}, and Proposition~\ref{prop-L-SL-has-FA}.
\end{proof}

\section{Boundaries and outer automorphisms}\label{sec-boundaries}
We prove the corollary for boundaries of random $\ell$-gonal groups.
\begin{proof}
	[Proof of Corollary 1.4]
	(1) As $G$ is a non-abelian free group, its boundary $\partial_\infty G$ is homeomorphic to a Cantor set.

	(2) As $G$ splits as $H \ast F_q$ for some $q \geq \frac{\sqrt{m}}{2}\geq 2$, its boundary has infinitely many connected components.  If its boundary was a Cantor set, then $G$ would be virtually free (see e.g.~\cite[Theorem 8.1]{kapovich-benakli-02-boundaries}).  But $G$ is torsion-free, so then $G$ would be free~\cite{Stallings68-torfree-ends}.  But, as in the proof of Theorem~\ref{thm-propFA-sharp}(2), the Euler characteristic of $G$ is a.a.s.\ $\chi(G) = 1-|S|+|R|>0$ and so $G$ is not free.  
	Thus $\partial_\infty G$ is neither connected nor a Cantor set.

	(3) Since every reduced diagram satisfies the isoperimetric inequality Corollary~\ref{cor-binomial-hyperbolicity}, the presentation complex for $G$ is aspherical, hence $G$ has cohomological dimension 2 (cf.\ \cite[Section I.3.b.]{Oll05}).  Thus $\partial_\infty G$ has topological dimension 1 \cite[Corollary 1.4(c)]{BestMess91}.  Since $G$ is torsion-free and has positive Euler characteristic  (by Lemma~\ref{lem-lots-of-relations}), by \cite[Remark at start of Section 6]{KapKle00} we have $\partial_\infty G$ is homeomorphic to a Menger sponge.
\end{proof}
Finally, we show the consequences for algebraic structure of $G$.
\begin{proof}
	[Proof of Corollary~\ref{cor-finite-out-sharp}]
	(1) 
	Let $\Aut(G)$ denote the automorphism group of $G$, and $\Inn(G) = \{ \alpha_g : g \in G\}$ the normal subgroup of inner automorphisms, $\alpha_g(h)=ghg^{-1}$, so $\Out(G) = \Aut(G)/\Inn(G)$.

	By Grushko's theorem (see e.g.~\cite{ScottWall79-top}), $G = G_1 \ast \cdots G_p \ast F_q$ for some freely indecomposable groups $G_1, \ldots, G_p$, $p \geq 0$, and some free group $F_q$, $q \geq 0$, and moreover any automorphism of $G$ preserves the $G_i$s up to conjugation and reordering.  
	In particular, if $N = \langle\langle G_1 \cup \cdots \cup G_p \rangle\rangle$, then any automorphism $\psi$ of $G$ preserves $N$, and we have a well-defined homomorphism $\pi:\Aut(G) \to \Aut(G/N)$.  
	If $\psi \in \Aut(G), g \in G$ then $\pi(\psi \alpha_g)(hN) = \pi(\psi)(ghNg^{-1}) = \pi(\psi)(ghg^{-1}N)$, so $\pi(\psi \alpha_g)$ and $\pi(\psi)$ represent the same element of $\Out(G/N) \cong \Out(F_q)$, and we have a well-defined homomorphism $\bar{\pi}:\Out(G) \to \Out(F_q)$.
	Moreover, $\bar{\pi}$ is surjective: given $[\phi] \in \Out(F_q)$, with $\phi \in \Aut(F_q)$, let $\hat{\phi}$ be the automorphism of $G$ given by letting $\hat{\phi}$ act as the identity on each $G_i$ and act as $\phi$ on $F_q$.  Then $\bar{\pi}([\hat{\phi}]) = \phi$.
	By Lemma~\ref{lem-missing-generators}, at this range of densities $q \geq \frac{\sqrt{m}}{2}\geq 2$.  
	Since $\Out(F_q)$ is infinite for $q \geq 2$, $\Out(G)$ is also infinite.

	We now show that $\mathrm{Hom}(G;\Gamma)$ is not finite up to conjugacy for any infinite torsion-free hyperbolic $\Gamma$: fix $x \in \Gamma$ a non-trivial element which generates a maximal infinite cyclic subgroup, and fix generators $y_1,\ldots,y_q$ of $F_q \leq G$.  Then for each $i \in \N$, let $\phi_i:G \to \Gamma$ be the homomorphism for which $\phi_i(y_1)=x^i$ and for which $\phi_i(y_j)=e$ for $j >1$, and $\phi_i(z) = e$ for all $z \in G_1 \cup \cdots \cup G_p$.  If $\phi_i = g \phi_j g^{-1}$ then $g$ lies in $E(x)$, the maximal elementary subgroup containing $x$ (compare \cite[Lemmas 1.16, 1.17]{Olshanskii-93-residualing-homomorphism-hyp-groups}).  By the choice of $x$ and the torsion-free assumption, $E(x) = \langle x\rangle$ so $g$ commutes with $x$ and we have $i=j$. Thus we have an infinite collection $\{\phi_i\}_{i \in \N}$ of elements of $\mathrm{Hom}(G;\Gamma)$ even after conjugation.

	(2) Since groups in this range have property FA by Theorem~\ref{thm-propFA-sharp}, this claim holds by the identical argument to \cite[Corollary 1.6]{DahGuiPrz11-nosplit}.  
\end{proof}

\bibliographystyle{alpha}
\bibliography{../biblio}

\begin{thebibliography}{BdlHV08}

\bibitem[A{\L}{\'S}15]{AntLucSwi-13threshold}
Sylwia Antoniuk, Tomasz {\L}uczak, and Jacek {\'S}wi{\c a}tkowski.
\newblock Random triangular groups at density 1/3.
\newblock {\em Compos. Math.}, 151(1):167--178, 2015.

\bibitem[ARD20]{AshcroftRonDou20}
Calum~J. Ashcroft and Colva~M. Roney-Dougal.
\newblock On random presentations with fixed relator length.
\newblock {\em Comm. Algebra}, 48(5):1904--1918, 2020.

\bibitem[Ash22]{Ashcroft22-14propT}
Calum~J Ashcroft.
\newblock Random groups do not have {P}roperty {(T)} at densities below 1/4.
\newblock arXiv preprint, 2022.
\newblock arXiv:2206.14616.

\bibitem[Ash23]{Ashcroft23-prop-T-lgonal}
Calum~James Ashcroft.
\newblock Property ({T}) in density-type models of random groups.
\newblock {\em Groups Geom. Dyn.}, 17(4):1325--1356, 2023.

\bibitem[BdlHV08]{BekkaHV09}
Bachir Bekka, Pierre de~la Harpe, and Alain Valette.
\newblock {\em Kazhdan’s Property (T)}.
\newblock New Mathematical Monographs. Cambridge University Press, 2008.

\bibitem[BM91]{BestMess91}
Mladen Bestvina and Geoffrey Mess.
\newblock The boundary of negatively curved groups.
\newblock {\em J. Amer. Math. Soc.}, 4(3):469--481, 1991.

\bibitem[DGP11]{DahGuiPrz11-nosplit}
Fran{\c c}ois Dahmani, Vincent Guirardel, and Piotr Przytycki.
\newblock Random groups do not split.
\newblock {\em Math. Ann.}, 349(3):657--673, 2011.

\bibitem[Gro96]{Gro93}
M.~Gromov.
\newblock Asymptotic invariants of infinite groups (geometric group theory).
\newblock {\em American Mathematical Society}, 33(3), July 1996.

\bibitem[J{\L}R00]{JansonLuczakRucinski-00-random-graphs-book}
Svante Janson, Tomasz {\L}uczak, and Andrzej Rucinski.
\newblock {\em Random graphs}.
\newblock Wiley-Interscience Series in Discrete Mathematics and Optimization.
  Wiley-Interscience, New York, 2000.

\bibitem[KB02]{kapovich-benakli-02-boundaries}
I.~Kapovich and N.~Benakli.
\newblock Boundaries of hyperbolic groups.
\newblock In {\em Combinatorial and geometric group theory (New York,
  2000/Hoboken, NJ, 2001)}, volume 296 of {\em Contemp.\ Math.}, pages 39--93.
  Amer.\ Math.\ Soc., Providence, RI, 2002.

\bibitem[KK00]{KapKle00}
Michael Kapovich and Bruce Kleiner.
\newblock Hyperbolic groups with low-dimensional boundary.
\newblock {\em Ann. Sci. \'Ecole Norm. Sup. (4)}, 33(5):647--669, 2000.

\bibitem[KK13]{Kotowski13}
Marcin Kotowski and Micha{\l} Kotowski.
\newblock Random groups and property {$(T)$}: {\.Z}uk's theorem revisited.
\newblock {\em J. Lond. Math. Soc. (2)}, 88(2):396--416, 2013.

\bibitem[MP15]{MacPrz15}
John~M. Mackay and Piotr Przytycki.
\newblock Balanced walls for random groups.
\newblock {\em Michigan Math. J.}, 64(2):397--419, 2015.

\bibitem[Odr16]{Odr16}
Tomasz Odrzyg\'o\'zd\'z.
\newblock The square model for random groups.
\newblock {\em Colloq. Math.}, 142(2):227--254, 2016.

\bibitem[Odr19]{Odr19}
Tomasz Odrzyg\'o\'zd\'z.
\newblock Bent walls for random groups in the square and hexagonal model.
\newblock arXiv preprint, 2019.
\newblock arXiv:1906.05417.

\bibitem[Ol'93]{Olshanskii-93-residualing-homomorphism-hyp-groups}
A.~Yu. Ol'shanski\u{i}.
\newblock On residualing homomorphisms and {$G$}-subgroups of hyperbolic
  groups.
\newblock {\em Internat. J. Algebra Comput.}, 3(4):365--409, 1993.

\bibitem[Oll05]{Oll05}
Yann Ollivier.
\newblock A {J}anuary 2005 invitation to random groups.
\newblock {\em Ensaios Matematicos [Mathematical Surveys], 10, Sociedade
  Brasileira de Mathemtica, Rio de Janeiro}, 2005.

\bibitem[OW11]{OllWise11}
Yann Ollivier and Daniel~T. Wise.
\newblock Cubulating random groups at density less than {$1/6$}.
\newblock {\em Trans. Amer. Math. Soc.}, 363(9):4701--4733, 2011.

\bibitem[Riv10]{Rivin10-growth-free-groups}
Igor Rivin.
\newblock Growth in free groups (and other stories)---twelve years later.
\newblock {\em Illinois J. Math.}, 54(1):327--370, 2010.

\bibitem[Ser03]{SerreTrees03}
Jean-Pierre Serre.
\newblock {\em Trees}.
\newblock Springer Monographs in Mathematics. Springer-Verlag, Berlin, 2003.
\newblock Translated from the French original by John Stillwell, Corrected 2nd
  printing of the 1980 English translation.

\bibitem[SPS85]{SchmidtPruzanShamir85}
Jeanette Schmidt-Pruzan and Eli Shamir.
\newblock Component structure in the evolution of random hypergraphs.
\newblock {\em Combinatorica}, 5(1):81--94, 1985.

\bibitem[Sta68]{Stallings68-torfree-ends}
John~R. Stallings.
\newblock On torsion-free groups with infinitely many ends.
\newblock {\em Ann. of Math. (2)}, 88:312--334, 1968.

\bibitem[SW79]{ScottWall79-top}
Peter Scott and Terry Wall.
\newblock Topological methods in group theory.
\newblock In {\em Homological group theory ({P}roc. {S}ympos., {D}urham,
  1977)}, volume~36 of {\em London Math. Soc. Lecture Note Ser.}, pages
  137--203. Cambridge Univ. Press, Cambridge-New York, 1979.

\bibitem[{\.Z}uk03]{Zuk03}
Andrzej {\.Z}uk.
\newblock Property ({T}) and {K}azhdan constants for discrete groups.
\newblock {\em (GAFA) Geometric and Functional Analysis}, 13(3):643--670, June
  2003.

\end{thebibliography}

\end{document}